\input amstex
\documentstyle{amsppt}
\input bull-ppt
\keyedby{bull406e/lbd}
\def\noi{{\noindent}}
\def\a{{\alpha}}
\def\b{{\beta}}
\def\th{{\theta}}
\def\d{{\delta}}

\topmatter
\cvol{29}
\cvolyear{1993}
\cmonth{July}
\cyear{1993}
\cvolno{1}
\shorttitle{Absence of Cantor spectrum}
\cpgs{85-87}
\title
Absence of Cantor spectrum for a class\\
of Schr\"odinger operators
\endtitle
\author
Norbert Riedel
\endauthor
\address
Department of Mathematics,
Tulane University,
New Orleans, Louisiana 70118\endaddress
\date July 15, 1992 and, in revised form, December 9, 
1992\enddate
\abstract
It is shown that the complete localization of
eigenvectors for the almost Mathieu operator
entails the absence of Cantor spectrum for this operator.
\endabstract
\endtopmatter

\document
\heading 1\endheading
\par
Among almost periodic Schr\"odinger operators the almost
Mathieu operator has enjoyed the majority of attention in 
the last ten to
fifteen years for the following three reasons: First, it 
is the simplest
nontrivial operator of its kind, and therefore it is 
deemed to be more
accessible than others.  Second, subjecting the almost 
Mathieu operator to
the Fourier transformation leads to an operator of the 
same kind (duality
property).  Third, and probably most important, depending 
on the parameters
involved, the almost Mathieu operator displays nearly all 
the conceivable
features operators of this kind could have, such as point 
spectrum, absence of
point spectrum, absolute continuous spectrum, singular 
spectrum, etc. According
to most researchers who have contributed to this field and 
who
have stated their opinion, however,
the only exception from this pattern seems to
concern the nature of the spectrum of the almost Mathieu 
operator, considered
as a bounded selfadjoint operator.  As manifested by 
repeatedly stated
conjectures (see [9] for one of the earliest sources and 
[5] for one of the
latest), the almost Mathieu operator is expected to have a 
nowhere dense
spectrum (Cantor spectrum) whenever the operator is not 
periodic.  We shall
argue that the absence of Cantor spectrum for the almost 
Mathieu operator not
only does occur, but we shall also identify the reason as 
to why it happens.
\bigskip
\heading 2\endheading
\par
The operator in question is defined on the Hilbert space
$\ell^2(\Bbb Z)$ as
$$(H(\alpha,\beta,\theta)\xi)_n=\xi_{n+1}+\xi_{n-1}+
2\beta\cos(2\pi\alpha
n+\theta)\xi_n\ ,$$
where $\a, \b, \th$ are real constants.  We will be 
concerned exclusively with
the case that $\a$ is an irrational number.  In this case 
the spectrum
$\roman{Sp}(\a,\b)$ of $H(\a,\b,\th)$ does not depend
on the parameter $\th$.  For a
full set of $\a$'s (in the sense of Lebesgue measure) 
which are sufficiently
badly approximable by rational numbers and for 
sufficiently large $\b$,
Fr\"ohlich et al. in [3] and Sinai in [10] have shown 
independently, with quite
%
involved techniques, that $H(\a,\b,\th)$ has, for almost 
all $\th$,  a
complete set of eigenvectors that decay exponentially.  
What is actually shown
in
these papers is that any generalized eigenvector of 
$H(\a,\b,\th)$ which grows
sufficiently slowly as $|n|$ approaches infinity decays 
exponentially. This
fact, combined with some other elementary observations, 
leads us to assume the
validity of the following condition for the parameters 
$\a$ and $\b$ in
question: \bigskip
\noi (L) \quad For every $\chi\in \roman{Sp}(\a,\b)$ there 
is a $\th$ such that
$H(\a,\b,\th)$ has an

\qquad eigenvector that decays exponentially as 
$|n|\to\infty$.
\bigskip
\noi This condition entails that the average Lyapunov 
index for the operator
$H(\alpha,\beta,\theta)$ takes the constant value 
$\log|\b|$ on $\roman{Sp}
(\a,\b)$.
By virtue of the ``showless formula", which establishes a 
connection between
the average Lyapunov index and the logarithmic potential 
associated with the
integrated density of states, it follows
that $\roman{Sp}(\a,\b)$ is a regular compactum (i.e., the 
Dirichlet problem
is solvable in $(\Bbb C \cup \{\infty\})\backslash 
\roman{Sp}(\a,\b)$
for any continuous
function on $\roman{Sp}(\a,\b)$), and the integrated 
density of states $\tau$ for
$H(\a,\b,\th)$ is the equilibrium distribution of 
$\roman{Sp}(\a,\b)$ (cf. [11] for
the  potential theoretic terminology being used).  The 
level curves of the
conductor potential associated with $\roman{Sp}(\a,\b)$ 
can be identified with the
spectra of perturbed almost Mathieu operators.
\proclaim {Theorem 1 \cite{7}}
Assume that {\rm(L)} holds.  A complex number $z$ is 
contained in the spectrum of the
operator $$(H_\d(\a,\b)\xi)_n=\xi_{n+1}+\xi_{n-1}+\b(\d 
e^{2\pi\a
ni}+\d^{-1}e^{-2\pi \a ni})\xi_n,\qquad \xi\in \ell^2(\Bbb 
Z)\ ,$$
if and only if $\int \log|z-s|\,d\tau(s)=\log|\b|+
|\log|\d||$.\endproclaim

The gist in the setup of the proof of this theorem is to 
consider the operators
$H_\d(\a,\b)$ as elements of the irrational rotation 
$C^*$-algebra ${\Cal
A}_\a$ associated with the number $\a$.  The $C^*$-algebra 
${\Cal A}_\a$ is
generated by two unitary operators $u$ and $v$ satisfying 
the commutation
relation $uv=e^{2\pi \a i}vu$.  The operator $H_\d(\a,\b)$ 
can be identified
with the element $u+u^{-1}+\b(\d v+\d^{-1}v^{-1})$.  The 
resolvent of the
operator $H_\d(\a,\b)$ on each connected component of the 
resolvent set can now
be expanded into absolutely convergent series in the 
monomials $u^pv^q$.
This in turn
leads to the definition of subharmonic functions measuring 
the order of
decay for these expansions in every point of the resolvent 
set.  The proof of
Theorem 1 is then established by exhibiting 
interrelationships between these
subharmonic functions. It should be noted that the 
conclusion of Theorem 1
remains true if the condition (L) is replaced by the 
assumption that $\a$ is
sufficiently well approximable by rationals and that 
$|\b|\ge 1$ \cite8. This
seems to suggest that the statement in Theorem 1 may be 
true for all irrational
numbers $\a$ and for all $|\b|\ge 1$.  By a basic duality 
argument it is
readily seen that a similar characterization of the level 
curves of the
conductor potential associated with $\roman{Sp}(\a,\b)$ in 
terms of the spectra of
perturbed almost Mathieu operators holds for any $\b$ with 
$0<|\b|<1$, whenever
the statement in Theorem 1 holds for $\b^{-1}$.
The condition (L), in conjunction with Theorem 1, is
indispensible for the following assertion.

\proclaim {Theorem 2 \cite7}  Assume
that {\rm(L)} holds.  Then $\roman{Sp}(\a,\b)$ is not a 
Cantor set.\endproclaim

The issue of the nature of $\roman{Sp}(\a,\b)$ has been 
taken up in a number of
papers.  Choi et al. [2], inspired by earlier work by 
Bellissard and Simon
[1], show that for numbers $\a$ which are sufficiently 
well approximable by
rationals, $\roman{Sp}(\a,\b)$ is indeed a Cantor set and,
moreover, they show that
all possible gaps in $\roman{Sp}(\a,\b)$ do actually 
occur.  The most penetrating
analysis to date of the nature of $\roman{Sp}(\a,\b)$ has 
been conducted by Helffer and
Sj\"ostrand in three memoirs [4].  Most of their 
investigation is limited to
the case $\b=1$ for $\a$'s which are sufficiently badly 
approximable by rational
numbers.  Even though this work produces a wealth of 
detailed information, it
leaves the main question unanswered.  In [10] Sinai claims 
that $\roman{Sp}(\a,\b)$ is
a Cantor set in those cases where he establishes complete 
localization of
eigenvectors.  As we have seen, this claim does not 
conform with our Theorem 2.

The proof of Theorem 2 heavily relies on some 
$C^*$-algebraic machinery that
was developed in [6] in order to address spectral problems 
for the almost
Mathieu operator.  In conjunction with Theorem 1, Theorem 2
provides for an argument
which translates the condition (L) into a smoothness 
condition for the
logarithmic potential involved.  In this context condition 
(L) appears as a
device to control the growth of the (noncommutative) 
Fourier expansions of
certain functionals associated with points in 
$\roman{Sp}(\a,\b)$.  That forementioned
smoothness of the logarithmic potential is shown to 
conflict with the
assumption that $\roman{Sp}(\a,\b)$ is a Cantor set.


\bigskip
\bigskip
\Refs
\rc
\ref
\no 1
\by J. Bellissard and B. Simon
\paper  Cantor spectrum for the almost Mathieu equation
\jour J. Funct. Anal.
\vol 48
\yr1982
\pages 408--519
\endref

\ref
\no 2
\by M.-D. Choi, G. A. Elliott, and N. Yui
\paper  Gauss polynomials and the rotation algebra
\jour Invent. Math.
\vol 99
\yr1990
\pages 225--246
\endref

\ref
\no 3
\by J. Fr\"ohlich, T. Spencer, and P. Wittwer
\paper Localization for a
class of one-dimensional quasi-periodic Schr\"odinger 
operators
\jour  Comm. Math. Phys.
\vol 132
\yr 1990
\pages 5--25
\endref

\ref
\no 4
\by B. Helffer and J. Sj\"ostrand
\book Analyse semi-classique pour l'\'equation de Harper. 
{\rm I-III}
\bookinfo M\'em. Soc. Math. France (N.S.)
(4) {\bf 116}, (4) {\bf 117}, (1) {\bf 118}\yr (1988--1990)
\endref

\ref
\no 5
\by L. Pastur and A. Figotin
\book Spectra of random and almost-periodic operators
\publ Springer-Verlag, New York
\yr 1992
\endref

\ref
\no 6
\by N. Riedel
\paper Almost Mathieu operators and rotation $C^*$-algebras
\jour  Proc. London Math. Soc. (3)
\vol 56
\yr 1988
\pages 281--302
\endref

\ref
\no 7
\bysame
\paperinfo {\it The spectrum of a class of almost periodic 
operators\/},
preprint
\endref

\ref
\no 8
\bysame
\paperinfo {\it Regularity of the spectrum for the almost 
Mathieu
operator\/}, preprint
\endref

\ref
\no 9
\by B. Simon
\paper Almost periodic Schr\"odinger
operators\,\RM: a review
\jour  Adv. in Appl. Math.
\vol  3
\yr 1982
\pages 463--490
\endref

\ref
\no 10
\by Ya. G. Sinai
\paper Anderson localization for one-dimensional
difference Schr\"odinger operator with quasiperiodic 
potential
\jour  J. Stat. Phys.
\vol 46
\yr 1987
\pages 861--909
\endref

\ref
\no 11
\by  M. Tsuji
\book Potential theory in modern function theory
\publ Chelsea, New York
\yr 1959
\endref
\endRefs
\enddocument